\documentclass[10pt]{amsart}
\pdfoutput=1

\usepackage{tikz}

\usepackage{amsmath}
\usepackage{url}
\usepackage{amssymb}
\usepackage{amsthm}
\usepackage{mathtools}
\usepackage{verbatim}
\usepackage{float}

\usepackage{caption}
\usepackage{graphicx}
\usepackage{graphics}

\usepackage{shortcuts}

\title[]{Many Lemniscates with large diameter}
\author[]{Linhang Huang}
\address{Department of Mathematics, University of Washington, Seattle}
\email{lhhuang@uw.edu}

\begin{document}
\begin{abstract}
We prove that for every $0 < c < 4$ and every $N \in \mathbb{N}$ there exists a monic polynomial $p(z) = z^n + a_{n-1} z^{n-1} + \dots + a_0$ such that the set $\set{z \in \C : |p(z)| \leq 1}$ has at least $N$ connected components with diameter at least $c$. This answers a question of Erd\H{o}s.
\end{abstract}
\maketitle

\section{Introduction and Result}

 Let $p(z) = z^n + a_{n-1} z^{n-1} + \dots + a_0$ be a monic polynomial in the complex plane. Henri Cartan (see, for example, \cite{levin}) proved that the set
 $ \left\{z \in \mathbb{C}: |p(z)| \leq 1\right\}$
can be covered by a union of circles whose sum of radii does not exceed $2e$. This has inspired many subsequent questions. In particular, Erd\H{o}s \cite{erdos} asked whether the set $ \left\{z \in \mathbb{C}: |p(z)| \leq 1\right\}$ can only contain few components whose diameter exceeds 1 or, more precisely, whether the number of connected components with diameter $>1+c$ is bounded from above by a universal constant $A(c)$ that is independent of the degree $n$. It is problem $\# 511$ in the collection of Erd\H{o}s problems \cite{bloom2}. 

\begin{center}
    \begin{figure}[h!]
\begin{tikzpicture}
    \node at (0,0) {\includegraphics[width=0.25\textwidth]{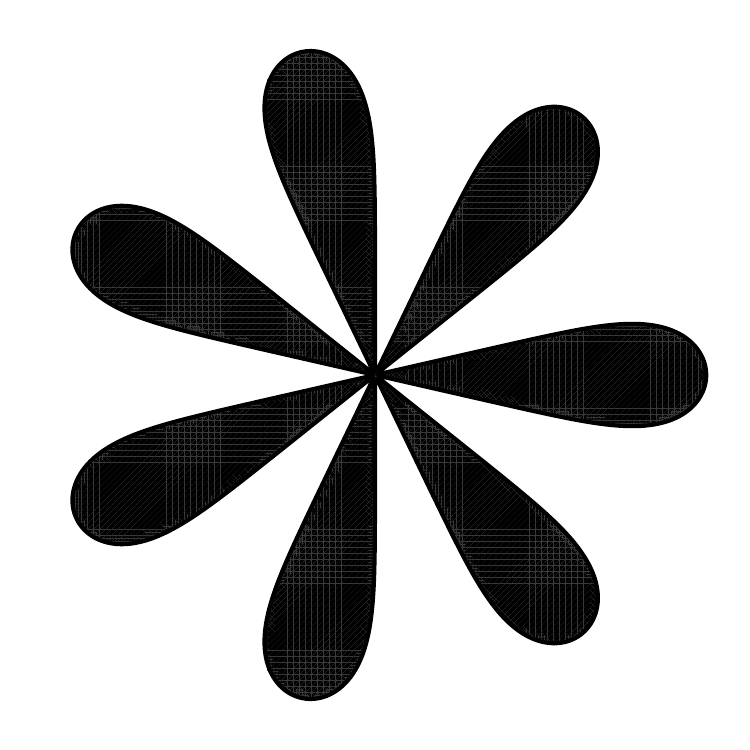}};
      \node at (6,0) {\includegraphics[width=0.25\textwidth]{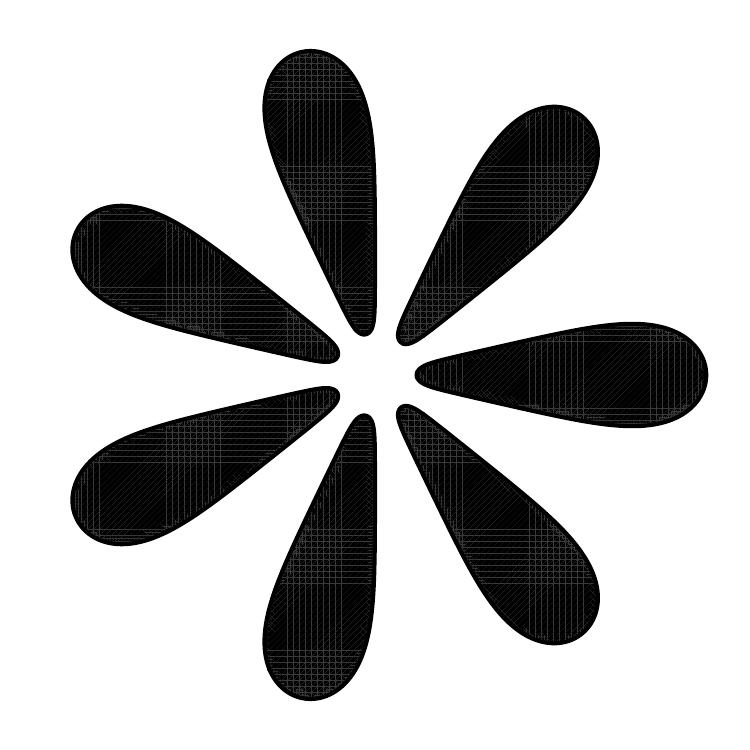}};
\end{tikzpicture}
\caption{The sets $|z^7 - 1| \leq 1$ and $|z^7 - 1| \leq 0.99999$.}
    \end{figure}
\end{center}
\vspace{-10pt}
One possible motivation is the polynomial $p(z) = z^n-1$ which is known to be extremal for a number of other problems. Then $|p(z)-1| \leq 1$ has one connected component, however, $|p(z) - 1| \leq 1-\varepsilon$
for $\varepsilon$ sufficiently small decouples into $n$ components of diameter slightly larger than 1 (depending on $n$).
We will show that this example is misleading: 
there exist examples where the diameter of 
each component is arbitrarily close to 4. 

\begin{theorem}
    For each $c\in (0,4)$ and $N \in \mathbb{N}$, there exists a monic polynomial
   $ p(z) = z^n + a_{n-1} z^{n-1} + \dots + a_0$
    such that $ \left\{z \in \mathbb{C}: |p(z)| \leq 1\right\}$ has at least $N$ connected components with diameter at least $c$.
\end{theorem} 

The restriction $c<4$ is best possible: P\'olya \cite{polya} showed that if $p:\mathbb{C} \rightarrow \mathbb{C}$ is a monic polynomial, then the orthogonal projection of $\left\{z \in \mathbb{C}: |p(z)| \leq 1\right\}$ onto any line is a set that can be covered by intervals whose length adds up to at most 4.
The main idea underlying our proof is to work with logarithmic capacity and to use the Hilbert Lemniscate Theorem to show that, for the purpose of this problem, there is no substantial difference between
$ \left\{z \in \mathbb{C}: |p(z)| \leq 1\right\}$ for a monic polynomial $p(z)$ and certain suitable sets with logarithmic capacity $1$. One direction is easy: the set
$ \left\{z \in \mathbb{C}: |p(z)| \leq 1\right\}$ has logarithmic capacity 1. The logarithmic capacity of a line segment of length $\ell$ is $\ell/4$ which is a way of suggesting $4$ as a fundamental limit on the diameter. There are sets with logarithmic capacity 1 whose diameter is arbitrarily close to 4 and these can be written down explicitly. We can then construct disjoint sets with diameter arbitrarily close to 1 whose union has capacity $\leq 1$. The Hilbert Lemniscate Theorem allows us to approximate such a set with the level set of a polynomial and capacity considerations imply that we may modify the polynomial to become monic.\\

\textbf{Note added.} When the problem was posted on \url{www.erdosproblems.com/511}, it was stated as unsolved. The author realized after submitting the initial preprint that the problem has been solved by Pommerenke in \cite{pommerenke}. This is an independent rediscovery. It will be kept on the arXiv for archival reasons but it will not be submitted to any journal.

\section{Proof}
The proof decouples into the following steps.
\begin{enumerate}
    \item We construct a family of domains with logarithmic capacity 1 and diameter $0 <c < 4$. The construction is fairly explicit.
    \item We build a union of $N$ disjoint Jordan curves with large diameter. There is a lot of freedom in this step.
    \item We apply the Hilbert Lemniscate Theorem to construct the polynomial and use capacity computations to show that it can be chosen to be monic.
\end{enumerate}

\subsection{A family of domains.}
For $0 < c < 4$ consider a shifted Joukowski map
$$ \varphi(z) = \frac{c}{4}\left( z + \frac{1}{z} + 2\right).$$

Note that the Joukowski map $z+1/z$ is a biholomorphic mapping between $\hat \C \backslash \overline\D$ and $\hat \C \backslash [-2,2]$. It follows that $\phi$ is a biholomorphic mapping between $\hat \C \backslash \overline\D$ and $\hat \C \backslash [0,c]$ and as $z \to \infty$, we have \begin{equation*}
    \varphi(z) = (c/4)z + \mathcal{O}(1).
\end{equation*}

\begin{figure}[h!]
   \begin{tikzpicture}[scale=1.3]
   \fill[domain=0:360,smooth,samples=30,variable=\t, gray, opacity=0.8]
  plot({cos(\t) + 2.5^2/16*cos(\t)+2.5/2},{sin(\t)-2.5^2/16*sin(\t)});

  \fill[domain=0:360,smooth,samples=30,variable=\t, gray, opacity=0.8]
  plot({cos(\t) + 3^2/16*cos(\t)+3/2+5},{sin(\t)-3^2/16*sin(\t)});

    \draw[->, line width = 1] (-0.5,0)--(3.5,0) node[right]{$x$};
    \draw[->, line width = 1] (0,-1.3)--(0,1.3) node[above]{$y$};
        \draw[red, line width = 2]
    (0,0) -- (2.5,0);
     \fill [black] (0,0) circle (1.2pt) node [black, below left=0.6pt]{$O$};

    \foreach \x in {1,2,3}
   \draw (\x cm,1pt) -- (\x cm,-1pt) node[anchor=north] {$\x$};
    \foreach \y in {-1,1}
    \draw (1pt,\y cm) -- (-1pt,\y cm) node[anchor=east] {$\y$};

    \draw[->, line width = 1] (4.5,0)--(8.5,0) node[right]{$x$};
    \draw[->, line width = 1] (5,-1.3)--(5,1.3) node[above]{$y$};
    
    \draw[red, line width = 2]
    (5,0) -- (8,0);
    
    \fill [black] (5,0) circle (1.2pt) node [black, below left=0.6pt]{$O$};

    \foreach \x in {1,2,3}
   \draw ({\x+5} ,1pt) -- ({\x+5} ,-1pt) node[anchor=north] {$\x$};
    \foreach \y in {-1,1}
    \draw (5cm+1pt,\y cm) -- (5cm-1pt,\y cm) node[anchor=east] {$\y$};

   \end{tikzpicture}
   \caption{The domains $\Omega$ (highlighted in gray) enclosed by $\varphi(\left\{|z| = 4/c\right\})$ for $c=2.5$ and $c=3$. As $c \rightarrow 4$, the curve approaches $[0,c]$ (highlighted in red).}
   \label{fig:twoc}
\end{figure}
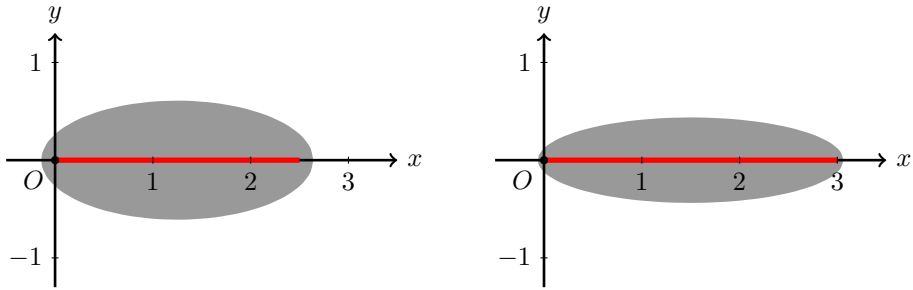

We will now consider the domain $\Omega$ enclosed by the curve $\varphi(\left\{|z| = 4/c\right\})$ which is completely explicit (see Fig. \ref{fig:twoc}). 
  $\Omega$ contains the interval $[0,c]$ and therefore has diameter at least $c$. The next step of the argument consists in showing that $\Omega$ has logarithmic capacity 1. Recall that for a compact set $K \subset \mathbb{C}$, the logarithmic capacity is defined as
    $$ \mbox{Cap}(K) = \exp\left( - \inf_{\mu} \int_{K \times K} \frac{1}{\log{|w-z|}} d\mu(w) d\mu(z) \right),$$
where the infimum ranges over all probability measures supported on $K$ (see Ransford \cite{ransford} for an introduction to logarithmic capacity). We will use an equivalent formulation in terms of the logarithmic capacity of connected compact sets in terms of biholomorphic mappings. More precisely, if $K$ is a compact set whose the unbounded component $U$ in $\hat \C$ is simply connected, then there exists a unique Riemann mapping function $g$ from $\hat\C \backslash \overline\D$ to $U$ such that $g(\infty) = \infty$ and $g'(\infty) >0$. And we have that \begin{equation*}
    \Cp(K) = g'(\infty)\quad \text{and thus}\quad g(z) = \Cp(K)z + \mathcal{O}(1).
\end{equation*}

 We now deduce $\Cp(\Omega) = 1$ from the fact that a rescaling of $\varphi$ satisfies
       $$\varphi(4z/c) = z + \mathcal{O}(1)$$
       while obviously being the Riemann mapping from $\hat\C \backslash \overline\D$ to $\hat{C} \backslash \overline{\Omega}$.  

\subsection{Constructing $\Omega_N$.}
\begin{figure}[h!]
   \begin{tikzpicture}[scale=0.6]
\fill[domain=0:360,smooth,samples=30,variable=\t, gray, opacity=0.8]
  plot({4.5*(cos(\t) + 2.2^2/16*cos(\t)+1.1)},{4.5*(sin(\t)-2.2^2/16*sin(\t))-0.3});
\draw[rounded corners = 4pt, black, line width=1.4] (-0.2, 1.1) rectangle (10,0.6);
\draw[rounded corners = 4pt, black, line width=1.4] (-0.2, 0.4) rectangle (10,-0.1);
\draw[rounded corners = 4pt, black, line width=1.4]  (-0.2, -0.5) rectangle (10,-1);
\draw[rounded corners = 4pt, black, line width=1.4]  (-0.2, -1.2) rectangle (10,-1.7);
\node at (6,-2.4) {\Huge $\partial \Omega_N$};

    \draw[red, line width = 2]
    (0,-0.3) -- ({2.2*4.5},-0.3);
     
   \end{tikzpicture}
   \caption{$\Omega$ has diameter $>c$ which allows us to find $\Omega_N$ enclosed by a union of $N$ disjoint Jordan curves of diameter $\geq c$.}\label{fig:jordan}
\end{figure}
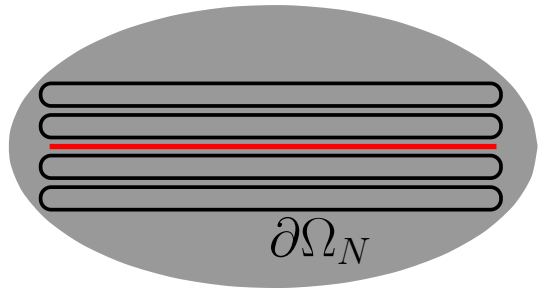

Let us now take $N$ pairwise disjoint Jordan curves inside $\Omega$ such that each enclosed domain has diameter greater than $c$. 
There are many ways one could do this and the precise way it is done does not matter: the subsequent steps in the argument only require that any two Jordan curves are separated by some positive quantity and that none of them touch the boundary of the domain. An example of what that could look like is shown in Fig. \ref{fig:jordan}. We denote the union of the corresponding $N$ domains bounded by the curves by $\Omega_N$. 

\subsection{The Hilbert Lemniscate Theorem.}

At this point, we recall the Hilbert Lemniscate Theorem \cite{hilbert}. A modern formulation, taken from  Bloom-Levenberg-Lyubarskii \cite{bloom}, is as follows: let $K \subset \mathbb{C}$ be a compact set with connected complement. For any $\varepsilon >0$ there exists a polynomial $p:\mathbb{C} \rightarrow \mathbb{C}$ such that
    $$K \subset \left\{ z \in \mathbb{C}: |p(z)| \leq \sup_{z \in K} |p(z)| \right\} \subset \left\{z \in \mathbb{C}: d(z, K) \leq \varepsilon\right\}.$$
We note that there is a (somewhat) constructive proof: one could take the polynomial to be the Fekete polynomial of sufficiently large degree corresponding to the set $K$ (see, for example, \cite{bloom}). We use the Hilbert Lemniscate Theorem with $K = \Omega_N$ and $\varepsilon$ so small that the
$\varepsilon-$neighborhood of $\Omega_N$ 
\begin{enumerate}
    \item has $N$ disjoint connected components that do not touch
    \item and is fully contained inside $\Omega$.
\end{enumerate}

Applying the Hilbert Lemniscate Theorem  with $K = \Omega_N$ and $\varepsilon$ guarantees the existence of a polynomial $q:\mathbb{C} \rightarrow \mathbb{C}$ such that 
\begin{align*}
     \Omega_N \subseteq \left\{z \in \mathbb{C}: |q(z)| \leq \sup_{z \in \Omega_N} |q(z)| \right\} \subseteq  \left\{z \in \mathbb{C}:  d(z,\Omega_N) \leq \varepsilon \right\} \subset \Omega.
\end{align*}

In particular, the set $\left\{z \in \mathbb{C}: |q(z)| \leq \sup_{z \in \Omega_N} |q(z)|\right\}$ has at least $N$ connected components with diameter at least $c$. After dividing the polynomial by $\sup_{z \in \Omega_N} |q(z)|$ and abusing notation by calling that polynomial again $q$, we have found a polynomial with the desired properties. However, this polynomial will not necessarily have leading coefficient 1 and it remains to show that we can choose the polynomial to be monic.
We use the following fact \cite[Theorem 5.2.5]{ransford}: if  
    $$ r(z) = a_d z^d + a_{d-1} z^{d-1} + \dots + a_0 \qquad \mbox{with}~a_d \neq 0,$$
    then
    $$ \mbox{Cap}(r^{-1}(\mathbb{D})) = |a_d|^{-1/d}.$$
    Since $q^{-1}(\overline{\D})\subseteq \Omega$, we deduce that $\Cp(q^{-1}(\overline{\D})) \leq \Cp(\Omega) = 1$. Therefore, the leading coefficient $a_d$ of $q$ satisfies $|a_d| \geq 1$. The desired result now follows from rescaling. Let $w \in \mathbb{C}$ be a solution of $w^d = a_d$ and consider the polynomial
   $ p(z) = q(z/w).$
    Then $p$ is monic. Moreover,
   \begin{equation*}
        p^{-1}(\overline\D) = w \cdot q^{-1}(\overline\D)
    \end{equation*} 
    which is simply a rotation and rescaling that does not change the topology. Since $|w| \geq 1$, we show that the diameters can only increase. This gives the desired result.


\begin{thebibliography}{10}

\bibitem{bloom} T. Bloom, N. Levenberg and Y. Lyubarskii, A Hilbert Lemniscate Theorem in $\mathbb {C}^ 2$. In Annales de l'institut Fourier 58 (2008), p. 2191-2220.

\bibitem{bloom1} T. F. Bloom, Erd\H{o}s Problem \#509, https://www.erdosproblems.com/509


\bibitem{bloom2} T. F. Bloom, Erd\H{o}s Problem \#511, https://www.erdosproblems.com/511

\bibitem{erdos} P. Erd\H{o}s: Some unsolved problems, Some unsolved problems. Magyar Tud. Akad. Mat. Kutat\'{o} Int. K\"{o}zl. (1961), 221-254. 

\bibitem{hilbert} D. Hilbert, \"Uber die Entwicklung einer beliebigen analytischen Funktion einer Variablen in eine unendliche nach ganzen rationalen Funktionen fortschreitende Reihe, Dritter Band: Analysis, Grundlagen der Mathematik, Physik, Verschiedenes: Nebst Einer Lebensgeschichte. Berlin, Heidelberg: Springer Berlin Heidelberg, 1935. 3-9.

\bibitem{levin} B. Ya. Levin, Distribution of zeros of entire functions, Transl. Math. Monographs, vol. 5, Amer. Math. Soc., Providence, RI, 1980.

\bibitem{polya} G. P\'olya, Beitrag zur Verallgemeinerung des Verzerrungssatzes auf mehrfach
zusammenhiingenden Gebieten, Sitzungsber. Preuss. Akad. Wiss. Berlin
(1928),228-232; Collected Papers Vol. I, MIT Press 1974, 347-351.

\bibitem{pommerenke} C. Pommerenke, On metric properties of complex polynomials, Michigan Math. J. {\bf 8} (1961), 97--115.

\bibitem{ransford} T. Ransford, Potential theory in the complex plane (No. 28). Cambridge University Press, 1995



\end{thebibliography}
\end{document}